\documentclass[12pt]{amsart}

\usepackage{amssymb}

\usepackage{enumerate}

\makeatletter
\@namedef{subjclassname@2010}{%
  \textup{2010} Mathematics Subject Classification}
\makeatother

\newcommand{\C}{\mathbb{C}}

\newcommand{\N}{\mathbb{N}}

\newtheorem{theorem}{Theorem}
\newtheorem{cor}{Corollary}
\newtheorem{lemma}{Lemma}

\newtheorem*{lemma*}{Lemma}
\newtheorem*{conj}{Conjecture}

\newtheorem{maintheorem}{Main Theorem}

\theoremstyle{definition}

\numberwithin{equation}{section}

\begin{document}


\baselineskip=17pt



\title[Averages of $\pi_{2k}(x)$]{Averaged Form of the Hardy-Littlewood Conjecture}

\author[J. Merikoski]{Jori Merikoski}
\email{jori.merikoski@helsinki.fi}

\date{\today}

\begin{abstract}
We study the prime pair counting functions $\pi_{2k}(x),$ and their averages over $2k.$ We show that good results can be achieved with relatively little effort by considering averages. We prove an asymptotic relation for longer averages of $\pi_{2k}(x)$   over $2k \leq x^\theta,$ $\theta > 7/12,$ and give an almost sharp lower bound for fairly short averages over $k \leq C \log x,$ $C >1/2.$ We generalize the ideas to other related  problems.
\end{abstract}

\subjclass[2010]{Primary 11P32 ; Secondary 11A41}

\keywords{Prime numbers, Hardy-Littlewood conjecture.}

\maketitle

\section{Introduction and Results}
In this article the main object of study is the counting function for prime number pairs
\begin{align}
\pi_{2k}(x)  := \vert \{ p \leq x: \, p \,\, \textrm{and} \,\, p+ 2k \,\, \textrm{both prime numbers} \} \vert,
\end{align}
where $2k \geq 2$ is an even integer. We use the notation $\vert A \vert$ for the cardinality of a given finite set $A.$  In particular, we discuss the asymptotic behaviour of  averages of the form
\begin{align*}
\frac{2}{M(x)} \sum_{2k \leq M(x)} \pi_{2k}(x),
\end{align*}
as $x$ tends to infinity. The structure of the article is as follows: In the first section we give a  brief introduction to the topic and state our main results. In the middle sections we give proofs for our main results. In the last sections we discuss generalizations, other related results, and conclusions. The ideas presented in this paper were conceived by the author during writing his Master's thesis at University of Helsinki. While the methods used are not very complicated, to our knowledge the results achieved do not exist in the literature.

Let us first recall the conjecture by Hardy and Littlewood \cite{HL} on the asymptotic behaviour of $\pi_{2k}(x).$
\begin{conj}\emph{\textbf{(Hardy-Littlewood Conjecture).}} \label{hl} Let $2k \geq 2$ be a constant. Then
\begin{align}
\pi_{2k}(x) \, \sim \, C_{2k} \frac{x}{\log^2 x}, 
\end{align}
as $x \to \infty,$ where the constant $C_{2k}$ is defined by
\begin{align}
C_{2k}: = 2 \prod_{p > 2} \frac{p(p-2)}{(p-1)^2} \prod_{2 < p \, \vert \, k} \frac{p-1}{p-2}.
\end{align}
\end{conj}

The conjecture remains open. In fact, it is still not known whether there are infinitely many prime pairs for any given $2k.$ The best result in this direction is by the recent online Polymath8 project, which states that for at least one even integer $2k \leq 246$ there are infinitely many primes $p$ such that $p+2k$ is also a prime number. That is,
\begin{align*}
\liminf_{k\to \infty} \, p_{k+1} - p_k \leq 246.
\end{align*}
The first result of this form was obtained by Yitang Zhang for at least one $2k \leq 70,000,000$ in 2013 \cite{Polymath}. It should be noted that the Prime Number Theorem $\pi(x) \sim x/\log x$ immediately implies that
\begin{align*}
\liminf_{k\to \infty} \, \frac{p_{k+1} - p_k}{\log p_k}\leq 1.
\end{align*}

Since the conjecture itself appears impregnable, we consider averages over $2k \leq M(x)$ for suitable functions $M(x) \to \infty$ as $x \to \infty.$ If the Hardy-Littlewood Conjecture holds uniformly for all $2k \leq M(x),$ then
\begin{align*}
\frac{2}{M(x)} \sum_{2k \leq M(x)} \pi_{2k}(x) \sim \frac{x}{\log^2 x } \frac{2}{M(x)} \sum_{2k \leq M(x)} C_{2k} \sim 2\frac{x}{\log^2x},
\end{align*}
as $x\to \infty$, where we have used the following lemma by Gallagher \cite{Gallagher}. It should be noted that Gallagher actually proved a more general version of the lemma.

\begin{lemma*} \emph{\textbf{(Gallagher).}} \label{galla} Let
\begin{align}
C_{2k} = 2 \prod_{p > 2} \frac{p(p-2)}{(p-1)^2} \prod_{2 < p \, \vert \, k} \frac{p-1}{p-2}
\end{align}
as in Conjecture \ref{hl}. Then
\begin{align}
\frac{2}{y} \sum_{2k \leq y} C_{2k} \to 2, \quad y \to \infty.
\end{align}
\end{lemma*}

The next theorem is our first main result. By the above we know that it is consistent with a uniform version of the Hardy-Littlewood Conjecture.

\begin{maintheorem}
\label{pairst}
Let $M(x)$ be a positive increasing function such that $M(x) \to \infty$ as $x \to \infty$ and $M(x) = o(x/\log^2 x).$ Suppose that for all functions $h= h(x),$ such that $M(x) \leq h(x) \leq x,$ we have \begin{align*}
\pi(x + h(x)) - \pi(x) \sim \frac{h(x)}{\log x}, \quad x \to \infty.
\end{align*}    Then 
\begin{align} \label{pairs}
\frac{2}{M(x)} \sum_{2k\leq M(x)} \pi_{2k}(x) \, \sim \, 2 \frac{x}{\log^2 x}, \quad x \to \infty.
\end{align}
Furthermore, for all such functions $h= h(x)$ we have
\begin{align}
\frac{2}{M(x)} \sum_{2k\leq M(x)} \left( \pi_{2k}(x +h) -   \pi_{2k}(x)\right) \, \sim \, 2 \frac{h}{\log^2 x}, \quad x \to \infty.
\end{align}
\end{maintheorem}

The function  $\pi(x + h(x)) - \pi(x)$ in the assumption of the previous theorem appears naturally when one tries to give upper bounds for the difference between consecutive prime numbers.   By Huxley \cite{Huxley}  we know that we may take $M(x) = x^{\theta} $ for any $\theta > 7/12$ in the above theorem, so that we have the following corollary.

\begin{cor}
Let $7/12 < \theta < 1.$  Then we have
\begin{align} \label{pairs}
\frac{2}{x^{\theta}} \sum_{2k\leq x^\theta} \pi_{2k}(x) \, \sim \, 2 \frac{x}{\log^2 x}, \quad x \to \infty.
\end{align}
Furthermore, for all functions $x^\theta \leq h \leq x$ we have
\begin{align}
\frac{2}{x^{\theta}} \sum_{2k\leq x^{\theta}} \left( \pi_{2k}(x +h) -   \pi_{2k}(x)\right) \, \sim \, 2 \frac{h}{\log^2 x}, \quad x \to \infty.
\end{align}
\end{cor}

It is worth noting about the main assumption in Theorem \ref{pairst} that if one is interested only in lower bounds of the form $ \pi(x + x^\theta) - \pi(x) \gg x^\theta/\log x,$ then we can choose even smaller $\theta = 0.525$  \cite{BHJ}. To obtain lower bounds for averages of $\pi_{2k}(x)$ over shorter intervals, we use a different method which does not depend on any such results.  To state our second main theorem, we need the following notation: For two real functions $f(x)$ and $g(x),$ we write $f(x) \geq g(x)(1+o(1))$ to say that for all $\epsilon > 0$ we have $f(x) \geq (1- \epsilon)g(x)$ for large enough $x \geq x_{\epsilon}.$ Similarly,  $f(x) \leq g(x)(1+o(1))$ means that  for all $\epsilon > 0$ we have $f(x) \leq (1+ \epsilon)g(x)$ for large enough $x \geq x_{\epsilon}.$ Our second main result gives a lower bound for a weighted average of $\pi_{2k}(x)$ over much shorter interval $2k \leq C \log x.$

\begin{maintheorem} 
Let $C > 1/2$ be a constant, and let $E =E(x)$ be such that $C \log x \leq E = o(x/\log^2x).$ Then
\begin{align*}
\frac{1}{\left\lfloor E \right\rfloor^2 }\sum_{1 \leq k \leq E} (\left\lfloor E \right\rfloor - k) \pi_{2k}(x) \geq \left(1-\frac{1}{2C}\right) \frac{ x}{\log^2 x}(1+o(1)).
\end{align*}
\end{maintheorem}

To see how this relates to the Hardy-Littlewood Conjecture, we compute using Abel's Summation Formula, and Gallagher's Lemma 
\begin{align*}
\frac{1}{E^2} \sum_{k\leq E} (E-k)C_{2k} & = \frac{1}{E^2}\int_1^E \sum_{k\leq x} C_{2k} \, dx \\
& = \frac{1}{E^2}\int_1^E (2x+o(x))\, dx \to 1, \quad E \to \infty.
\end{align*}
Hence, for $E$ such that $\log x = o(E),$ the lower bound in the above theorem is sharp assuming that the Hardy-Littlewood Conjecture holds uniformly for $2k \leq E.$

As an immediate corollary we get the following  more elegant but weaker version.
\begin{cor}
Let $C$ and $E $ be as in Theorem 2. Then
\begin{align*}
\frac{1}{\left\lfloor E \right\rfloor} \sum_{1 \leq k \leq E}  \pi_{2k}(x) \geq \left(1-\frac{1}{2C}\right) \frac{ x}{\log^2 x}(1+o(1)).
\end{align*}
\end{cor}

Since $E$ can be taken close to $\log x /2,$ this corollary is a much more quantitative statement of the fact that  
\begin{align*}
\liminf_{k\to \infty} \, \frac{p_{k+1} - p_k}{\log p_k}\leq 1.
\end{align*} 
As mentioned before, this follows immediately from the Prime Number Theorem. 

\section{Proof of Main Theorem 1}
Let $P(n)$ denote the characteristic function of primes. That is,
\begin{align}
P(n) := \begin{cases} 1, \, & n \,\, \textrm{ is a prime,} \\
0 \, & \textrm{otherwise.}
\end{cases}
\end{align} 
Then for any fixed $2k$ the function $P(n)P(n+2k)$ is the characteristic function for primes $p$ such that $p+2k$ is also a prime. Note that if $r$ is odd, then the number of primes $p$ such that $p+r$ is also a prime is at most $1.$ Therefore,
\begin{align*}
\sum_{2k\leq M(x)} \pi_{2k}(x) & =  \sum_{2k\leq M} \sum_{n \leq x} P(n)P(n+2k) \\
& =  \sum_{k\leq M} \sum_{n \leq x} P(n)P(n+k) + \mathcal{O}(M) \\
& = \sum_{k\leq M} \sum_{ M < n \leq x} P(n)P(n+k) + \mathcal{O}(M^2).
\end{align*}
 Changing the order of summation yields
\begin{align*}
 \sum_{k\leq M} \sum_{ M < n \leq x} P(n)P(n+k) & =  \sum_{M< n \leq x} P(n) \sum_{k\leq M}  P(n+k) \\
& =  \sum_{ M < n \leq x} P(n) \left( \pi(n + M) - \pi(n) \right).
\end{align*}
Since $M< n \leq x$, we have $M(n) \leq M(x) < n,$ since $M$ is increasing. Hence, by our main assumption we have $\pi(n + M(x)) - \pi(n) \sim M(x)/ \log n,$ as $n,x \to \infty.$
Hence, the last sum is
\begin{align*}
 \sum_{ M  < n \leq x} P(n) \left( \frac{M(x)}{\log n} + o\left(  \frac{M(x)}{\log n} \right) \right) &= \sum_{ M< n \leq x} P(n)  \frac{M(x)}{\log n} + o\left( \max \{ M, \sqrt{x}\} +  \sum_{ \sqrt{x}  < n \leq x}  P(n) \frac{M(x)}{\log n} \right) \\
 &= \sum_{ M< n \leq x} P(n)  \frac{M(x)}{\log n} + o\left(  \frac{x M(x)}{\log^2 x} \right) 
\end{align*}
 by the Prime Number Theorem. The above sum is
\begin{align*}
\sum_{ M< n \leq x} P(n) \frac{M}{\log n} \geq \frac{M}{\log x} \left( \pi(x) -  \pi(M)  \right) \sim \frac{xM}{\log^2 x}.
\end{align*}
To obtain an inequality to the other direction let $M< y< x.$ Then
\begin{align*}
\sum_{ M < n \leq x} P(n)  \frac{M}{\log n}& \leq  \sum_{ y < n \leq x} P(n) \frac{M}{\log y} + \sum_{M< n \leq y} P(n) \frac{M}{ \log n}  \\
&  = \frac{M}{ \log y} \left( \pi(x) -  \pi(y)  \right) +  \sum_{\sqrt{x} < n \leq y} P(n) \frac{M}{ \log n}   + \mathcal{O} (\max \{ M^2, M\sqrt{x}\} ) \\
&  = \frac{M}{ \log y} \left( \pi(x) -  \pi(y)  \right) + \mathcal{O} \left(  \frac{y M}{\log x \log y} +\max \{ M^2, M\sqrt{x}\} \right) \\
& = \frac{xM}{\log y \log x} + \mathcal{O}\left(  \frac{My}{\log x} + \max \{ M^2, M\sqrt{x}\}\right)  \\
& =  \, \frac{xM}{\log^2 x} +   o \left(  \frac{ xM}{\log^2 x} \right)
\end{align*}
for $y=x/\log^2 x,$ since $M = o(x/\log^2 x).$
Hence,
\begin{align*}
\sum_{2k\leq M} \pi_{2k}(x) & =  \sum_{ M< n \leq x} P(n) \left( \frac{M}{\log n}\right) + o\left(  \frac{xM}{\log^2 x} \right)  + \mathcal{O}(M^2) \\
 & =\, \frac{xM}{\log^2 x}  +  o \left(  \frac{ xM}{\log^2 x} \right) + \mathcal{O}(M^2).
\end{align*}
Since $M(x) = o(x / \log^2 x),$ we obtain
\begin{align*}
\frac{2}{M} \sum_{2k\leq M} \pi_{2k}(x) \, \sim \, 2 \frac{x}{\log^2 x}
\end{align*}
as $ x \to \infty.$

The proof of the second claim in Main  Theorem 1 is similar, but slightly easier. For all $M(x) \leq h \leq x$ we have
\begin{align*}
\frac{2}{M} \sum_{2k\leq M} \left( \pi_{2k}(x +h) -   \pi_{2k}(x)\right)  \, &= \, \frac{2}{M} \sum_{2k\leq M} \sum_{x < n \leq x + h} P(n)P(n+2k) \\
& = \,  \frac{2}{M} \sum_{x < n \leq x + h} P(n) \sum_{2k\leq M}  P(n+2k) \\
& = \,  \frac{2}{M} \sum_{x < n \leq x + h} P(n) \left(  \pi(n + M) - \pi(n) \right) \\
& = \,  \frac{2}{M} \sum_{x < n \leq x + h} P(n)\left(\frac{M}{\log n} + o\left(\frac{M}{\log n} \right) \right) \\
& = \frac{2}{M} \left(\frac{h}{\log x} + o \left(\frac{h}{\log x} \right) \right)\left(\frac{M}{\log x} + o\left(\frac{M}{\log x} \right) \right)   \\
&\sim 2 \frac{h}{\log^2 x},
\end{align*}
where we have used the main assumption twice, and the fact that $\log n \sim \log x$ for all $x< n \leq x+h \leq 2x.$ \qed

\section{Proof of Main Theorem 2}
To prove the second main theorem, we need the following lemma which follows from a simple application of the Cauchy-Schwarz inequality and the Prime Number Theorem

\begin{lemma} \label{key}
Let $B = B(x)$ be a set of positive integers such that for all $b \in B$ we have $b = o\left( x/\log^2 x \right)$ uniformly. Then
\begin{align}
\frac{1}{\vert B \vert^2} \sum_{(a,b) \in B^2, \, a \neq b} \pi_{\vert a-b \vert}(x) \geq \left(\frac{x}{\log^2 x} - \frac{1}{\vert B \vert}\frac{x}{\log x} \right)(1+o(1)).
\end{align}
\end{lemma}
\begin{proof}
We have
\begin{align*}
\sum _{a \in B } \sum_{n \leq x} P(n+a) = \sum _{a \in B } (\pi(x+a) - \pi(a)) = \vert B \vert \pi(x) (1+o(1)),
\end{align*}
since $a = o\left(x/\log^2 x \right) = o\left( \pi(x) \right)$ for all $a \in B.$
Therefore, by the Cauchy-Schwartz inequality
\begin{align*}
\vert B \vert^2 \pi(x)^2 (1+o(1)) = \left( \sum_{n \leq x} \sum _{a \in B } P(n+a)  \right)^2 \leq  x \sum_{n \leq x}\left(  \sum _{a \in B } P(n+a)  \right)^2.
\end{align*}
Using the Prime Number Theorem this implies
\begin{align} \label{eqq}
\frac{1}{\vert B \vert^2}\sum_{n \leq x}  \sum_{(a,b) \in B } P(n+a) P(n+b)   \geq \frac{x}{\log^2 x}(1+o(1)).
\end{align}
The left-hand side is
\begin{align*}
\frac{1}{\vert B \vert^2}\sum_{(a,b) \in B^2} \left( \pi_{\vert a-b \vert}(x + \min \{a,b\}) - \pi_{\vert a-b \vert}(\min \{a,b\})\right),
\end{align*}
which is equal to
\begin{align*} 
 \frac{1}{\vert B \vert^2} \sum_{(a,b) \in B^2} \pi_{\vert a-b \vert}(x) + o \left( \frac{x}{\log^2 x}\right),
\end{align*}
since $a,b = o\left(x/\log^2 x\right).$
The contribution from the pairs $(a,b)$ such that $a = b$ is
\begin{align*}
\frac{1}{\vert B \vert} \pi(x) = \frac{1}{\vert B \vert} \frac{x}{\log x}(1+o(1)).
\end{align*}
Moving this to the right-hand side of \eqref{eqq} yields
\begin{align*}
\frac{1}{\vert B \vert^2} \sum_{(a,b) \in B^2, \, a \neq b} \pi_{\vert a-b \vert}(x) \geq \left(\frac{x}{\log^2 x} - \frac{1}{\vert B \vert}\frac{x}{\log x} \right)(1+o(1)).
\end{align*}
\end{proof}

We can now prove Theorem 2 by applying the above lemma to the set $B = \{ 1,2,3, \dots,  2 \lfloor E \rfloor \}.$ This yields
\begin{align*}
\frac{1}{4\lfloor E \rfloor ^2 } \sum_{(a,b) \in B^2, \, a \neq b} \pi_{\vert a-b \vert}(x) \geq \left(\frac{x}{\log^2 x} - \frac{1}{2C}\frac{x}{\log^2 x} \right)(1+o(1)),
\end{align*}
since $E \geq C \log x.$

If $a-b$ is not divisible by 2, then $\pi_{\vert a-b\vert}(x) \leq 1.$ The contribution from this to the sum on the left-hand side is clearly negligible.

Every even number $2k \leq 2\left\lfloor E \right\rfloor$ appears $4(\left\lfloor E \right\rfloor - k)$ times as the difference $\vert a-b \vert,$ namely for the pairs
\begin{align*}
(1, \, 1 + 2k), \, (2, \, 2+2k), \dots, \, (2\lfloor E \rfloor) - 2k, \, 2\lfloor E \rfloor),
\end{align*}
and  the other way around.
 Hence,
\begin{align*}
\frac{1}{\lfloor E \rfloor^2}\sum_{1 \leq k \leq E } (\lfloor E \rfloor - k) \pi_{2k}(x) \geq   \left(1-\frac{1}{2C} \right)\frac{ x}{\log^2 x}(1+o(1)).
\end{align*} \qed

\section{Generalizations and Related Results}
In this section we skecth some possible ways to generalize our results without  discussing too much  about technical details. One clear way to generalize our ideas is to consider other prime constellations than prime pairs. For example, for a sequence of even integers $(2h_1, 2h_2, \dots, 2h_k )$ we can define \begin{align*}
\pi_{2h_1, \dots, 2h_k}(x) =  \vert \{ p \leq x: \, p, p+2h_1, \dots, p+2h_k \, \, \textrm{all prime numbers} \} \vert.
\end{align*}  
We can  then study the averages of $\pi_{2h_1, \dots, 2h_k}(x)$ over $2h_1, \dots, 2h_k \leq M(x).$ Following the lines of the proof of Main Theorem \ref{pairst} we would obtain
\begin{align*}
\frac{2^k}{M^k} \sum_{2h_1\leq M}\cdots\sum_{2h_k\leq M}  \pi_{2h_1, \dots, 2h_k}(x) \,  & \sim \, \frac{2^k}{M^k}  \sum_{ M < n \leq x} P(n) \sum_{h_1 \leq M}  P(n+h_1) \cdots\sum_{h_k \leq M}  P(n+h_k) \\
& \sim \, 2^k \frac{x}{\log^{k+1} x},
\end{align*}
given the assumption that for all functions $h$  such that $M(x) \leq h(x) \leq x$ we have \begin{align*}
\pi(x + h) - \pi(x) \sim \frac{h}{\log x}, \quad x \to \infty.
\end{align*}  

The Hardy-Littlewood Conjecture generalizes to  $\pi_{2h_1, \dots, 2h_k}(x)$ as follows:  Let \begin{align*}
C_{2h_1, \dots, 2h_k } := \prod_{p }  \left( 1- \frac{\nu_{2h_1, \dots, 2h_k}(p)}{p} \right) \left( 1-\frac{1}{p}\right)^{-k-1},
\end{align*}
where $\nu_{2h_1, \dots, 2h_k}(p)$ is the size of the set $\{0, 2h_1, \dots, 2h_k\}$ modulo $p,$ that is, the number of residue classes of $p$ that the set $\{0, 2h_1, \dots, 2h_k\}$ meets. Note that for $k=1$ this agrees with our earlier definition. If for all prime numbers $p$ the set $\{0, 2h_1, \dots, 2h_k\}$ avoids at least one of the residue classes modulo $p,$   then we expect that
\begin{align*}
\pi_{2h_1, \dots, 2h_k}(x) \sim \,  C_{2h_1, \dots, 2h_k } \frac{x}{\log^{k+1} x}.
\end{align*}
This is known as the $\emph{Hardy-Littlewood k-tuple Conjecture}.$
Notice that if $\{0, 2h_1, \dots, 2h_k\}$ contains at least one member of each residue class modulo some $p,$ then  $C_{2h_1, \dots, 2h_k } =0,$ and $\pi_{2h_1, \dots, 2h_k}(x) $ is bounded. 

The general version of Gallagher's Lemma (see \cite{Gallagher}) implies that
\begin{align*}
\frac{2^k}{y^k} \sum_{2h_1\leq y}\cdots\sum_{2h_k\leq y} C_{2h_1, \dots, 2h_k } \to 2^k,
\end{align*}
as $y \to \infty.$ Hence, our generalization of Main Theorem 1 is consistent with a uniform version of the Hardy-Littlewood $k$-tuple Conjecture.

The method used to prove Main Theorem 2 also generalizes for prime $k$-tuples. In Lemma \ref{key}, we just need to replace Cauchy-Schwarz inequality  by H\"older's inequality. Repeating the argument, for a set $B$ such that $b = o(x/\log^{k+1}x)$ for all $b \in B$, we have
\begin{align*}
\sum _{a \in B } \sum_{n \leq x} P(n+a) = \sum _{a \in B } (\pi(x+a) - \pi(a)) = \vert B \vert \pi(x) (1+o(1)).
\end{align*}
Using H\"older's inequality we obtain
\begin{align*}
\vert B \vert^{k+1} \pi(x)^{k+1} (1+o(1)) &= \left( \sum_{n \leq x} \sum _{a \in B } P(n+a)  \right)^{k+1}  \\
& \leq  x^k \sum_{n \leq x}\left(  \sum _{a \in B } P(n+a)  \right)^{k+1} \\
& = x^k \sum_{n \leq x} \, \, \sum _{a_1, \dots, a_{k+1} \in B } P(n+a_1) \cdots P(n+a_{k+1}) \\
& = x^k  \sum _{a_1, \dots, a_{k+1} \in B } \pi_{ (a_1 - m), \dots, (a_{k+1} - m)} (x) + o\left(\frac{x^{k+1}}{\log^{k+1} x} \right),
\end{align*}
where $m := \min \{a_1, a_2, \dots a_{k+1} \}.$ The contribution of the terms where $a_1 = a_2 \cdots = a_{k+1}$ to the sum on the right-hand side is $x^k \vert B \vert \pi(x).$ Moving this to the other side, and using the Prime Number Theorem we obtain
\begin{align*}
\frac{1}{\vert B \vert^{k+1}}\sum_{\substack{a_1, \dots, a_{k+1} \in B \\ \exists i,j:\,  a_i \neq a_j }} \pi_{(a_1 - m), \dots, ( a_{k+1} - m)} (x) \geq \left( \frac{x}{\log^{k+1} x} - \frac{1}{\vert B \vert^k}\frac{x}{\log x}\right)(1+o(1)).
\end{align*}
Suppose then that $C > \frac{1}{2}$ and $E = E(x) \geq C \log x.$ Applying the above inequality to the set $B= \{1,2,3, \dots, 2\lfloor E \rfloor\}$ yields 
\begin{theorem}\begin{align*}
\frac{k+1}{2^k \lfloor E \rfloor^{k+1}} \sum_{h_1\leq E}\cdots\sum_{h_k\leq E} (\lfloor E \rfloor - M) \pi_{2h_1, \dots, 2h_k}(x)  \geq \left( 1- \frac{1}{(2C)^{k}}\right)\frac{x}{\log^{k+1}x}(1+o(1)),
\end{align*}
where $M= M(h_1, \dots, h_{k}) = \max \{h_1, \dots, h_k \}.$
\end{theorem}
\begin{proof}
In the sum \begin{align*}
\sum_{\substack{a_1, \dots, a_{k+1} \in B \\ \exists i,j:\,  a_i \neq a_j }} \pi_{(a_1 - m), \dots, ( a_{k+1} - m)} (x)
\end{align*}
the sequence $2h_1, \dots, 2h_k$ appears $2(k+1)(\lfloor E \rfloor - M)$ times, namely when $(a_1, \dots, a_{k+1})$ equals to
\begin{align*}
&(1,1+2h_1, \dots, 1+ 2h_k), (2,2+2h_1, \dots, 2+ 2h_k) , \dots \\ & \dots, (2\lfloor E \rfloor - 2M, 2\lfloor E \rfloor - 2M + 2h_1, \dots, 2\lfloor E \rfloor - 2M + 2h_k),
\end{align*}
if we fix $m = \min \{a_1, a_2, \dots a_{k+1} \}= a_1,$ and similarly for the $k$ other possible choices of $m.$ This gives the constant $2(k+1)(\lfloor E \rfloor - M).$
\end{proof}

For the record we note that Lemma \ref{key} holds in general for other arithmetic functions besides $P(n).$ Let $A: \N \to \C$ be a function and define $\alpha(x) := \sum_{n\leq x} A(n),$ $\alpha_k (x) := \sum_{n \leq x} A(n)\overline{ A(n+k)},$ where $k$ is an integer. Then the same argument as in Lemma \ref{key} yields
\begin{lemma} 
Let $B = B(x)$ be a set of positive integers such that for all $b \in B$ we have $b = o\left( \vert \alpha (x) \vert ^2 /x \right)$ uniformly. Then
\begin{align*}
\frac{1}{\vert B \vert^2} \left \vert \sum_{(a,b) \in B^2, \, a \neq b} \alpha_{a-b}(x)  \right \vert \geq \left(\frac{\vert \alpha (x) \vert ^2 }{ x} - \frac{\alpha_0(x)}{\vert B \vert}\right)(1+o(1)).
\end{align*}
\end{lemma} 
The above lemma is non-trivial only for sets $B$ such that  $\vert B \vert > x \alpha_0(x) / \vert \alpha (x) \vert ^2 .$ This can  be generalized further by using H\"older's inequality in a similar fashion as above, but we do not pursue this any further.

To conclude this section, we give an example of how Lemma \ref{key} can be used with a different kind of set $B.$ As a result,  we  obtain a lower bound for the of average of $\pi_{2mk}(x)$ over $k.$
\begin{theorem}
Let $h=h(x) = o\left(x/\log^2 x \right)$ such that $\log x = o(h).$ Let $m = o \left( h/\log x\right)$ be an integer. Then
\begin{align*}
\frac{1}{M^2} \sum_{ 1 \leq k \leq M} 2 \left(  M - k \right) \pi_{2mk}(x) \geq \frac{x}{\log^2 x}(1+o(1)),
\end{align*}
where $M =  \left\lfloor h/(2m) \right \rfloor.$
\end{theorem}
\begin{proof}
Let $B := \{ 2m, 4m, \dots, \, 2m M \}.$ Then by Lemma \ref{key} 
\begin{align*}
\frac{1}{M^2}\sum_{(a,b) \in B^2, \, a \neq b} \pi_{\vert a-b \vert}(x) \geq \left(\frac{x}{\log^2 x} - \frac{1}{M}\frac{x}{\log x} \right)(1+o(1)) = \frac{x}{\log^2 x}(1+o(1)),
\end{align*}
since $1/M \sim \, 2m/h = o \left( 1/\log x\right).$ Each number $2mk$ appears $2(M-k)$ times as the difference $\vert a-b \vert,$ which proves the theorem.
\end{proof} 
 The following more appealing version follows at once.
\begin{cor} 
Let $h=h(x)=o\left(x/\log^2 x \right)$ such that $\log x = o(h).$ Let $m = o \left( h/\log x \right)$ be a positive integer. Then
\begin{align*}
\frac{1}{M} \sum_{ 1 \leq k \leq M}  \pi_{2mk}(x) \geq \frac{x}{2\log^2 x}(1+o(1)),
\end{align*}
where $M =  \left\lfloor h/(2m) \right \rfloor.$
\end{cor}

\section{Conclusions}
We have studied the averages of $\pi_{2k}(x)$ over $2k,$ and shown that good results can be obtained with relatively little effort. We have also generalized these ideas to prime $k$-tuples. In particular, we have shown that long averages satisfy
\begin{align*}
\frac{2}{x^\theta} \sum_{2k\leq x^\theta} \pi_{2k}(x) \, \sim \, 2 \frac{x}{\log^2 x}, \quad x \to \infty
\end{align*}
for any $7/12 < \theta < 1$ which is consistent with the Hardy-Littlewood Conjecture by Gallagher's Lemma. 

In addition, we have shown that averages over fairly short intervals satisfy lower bounds
\begin{align*}
\frac{1}{\left\lfloor E \right\rfloor^2 }\sum_{1 \leq k \leq E} (\left\lfloor E \right\rfloor - k) \pi_{2k}(x) \geq \left(1-\frac{1}{2C}\right) \frac{ x}{\log^2 x}(1+o(1)),
\end{align*}
for $E \geq C \log x,$ $C \geq 1/2.$ For $E$ such that $\log x = o(E)$ this implies
 that
 \begin{align*}
\frac{1}{\left\lfloor E \right\rfloor^2 }\sum_{1 \leq k \leq E} (\left\lfloor E \right\rfloor - k) \pi_{2k}(x) \geq \frac{ x}{\log^2 x}(1+o(1)),
\end{align*}
which is the best possible bound if the Hardy-Littlewood conjecture holds uniformly. A topic for future investigations would be to obtain a similar upper bound for short averages, for example using sieve methods, which have been very effective in obtaining upper bounds (see e.g. \cite{HR}). This could be used together with our results to prove the following likely conjecture.
\begin{conj} For functions $E=E(x)$ such that $\log x = o(E)$ we have 
\begin{align*}
\frac{1}{\left\lfloor E \right\rfloor^2 }\sum_{1 \leq k \leq E} (\left\lfloor E \right\rfloor - k) \pi_{2k}(x) \sim \frac{ x}{\log^2 x}, \quad x \to \infty.
\end{align*}
\end{conj}

Another possible topic for future papers would be to use sieve theory or other methods to extend our results for shorter intervals.

\end{document}